\documentclass[11pt]{article}
\usepackage{amsfonts}
\usepackage{eufrak}
\usepackage{mathrsfs}
\usepackage{amssymb}
\usepackage{mathtools}
\usepackage{commath}
\usepackage{dsfont}
\begin{document}
\title{\bf On sums of unequal powers of primes\\ and powers of 2
\thanks{Project supported by the National Natural Science Foundation of China (grant no.11771333) and by the Fundamental Research Funds for the Central Universities (grant no. JUSRP122031). }}
\author{{Yuhui  Liu}\\ {School of Science,  Jiangnan University}\\{Wuxi 214122, Jiangsu, China}
\\{Email: tjliuyuhui@outlook.com}}
\date{}
\maketitle
{\bf Abstract}
In this paper, it is proved that every sufficiently large even integer can be represented as the sum of two squares of primes, two cubes of primes, two biquadrates of primes and 16 powers of 2. Furthermore, there are at least 5.313\% odd integers that can be represented as one square of prime, one cube of prime and
one biquadrate of prime. This result constitutes a refinement upon that of R. Zhang [8].



\section{Introduction}
\setcounter{equation}{0}
\hspace{5.8mm}In the 1950s, Linnik [2,3] proved that each large even integer $N$ is a sum of two primes and a bounded number of powers of 2,
\begin{align}
N = p_1 + p_2+ 2^{v_1} +2^{v_2} + \cdots + 2^{v_{k_1}},
\end{align}
where and below the letter $p$ and $v$, with or without subscripts, denote a prime number and a positive integer respectively. The famous Goldbach conjecture implies that $k_1=0$. The explicit value for the number $k_1$ was improved by many authors.

In 2017, Liu [5] considered a Goldbach-Linnik problem with unequal powers of primes, i.e.
\begin{align*}
N = p_1^2 + p_2^2 + p_3^3 + p_4^3 + p_5^4 + p_6^4 + 2^{v_1} + 2^{v_2} + \cdots + 2^{v_{k}},
\end{align*}
and proved that every sufficiently large even integer can be written as a sum of two squares of primes, two cubes of primes, two fourth powers of primes and at most 41 powers of 2. In 2019, L\"{u} [6] improved the value  to 24. Recently, the result was refined to 22 and then to 20, by Zhao [10] and Zhang [8], respectively. In this paper, we obtain a further improvement of  the value of $v_k$ by giving the following theorem.

{\bf Theorem 1.1} {\it \,\,Every sufficiently large even integer is a sum of two squares of primes, two cubes of primes, two biquadrates of primes and $16$ powers of $2$.}\\

The second result in this paper which involves unequal powers of primes is given in the positive density form.

{\bf Theorem 1.2} {\it \,\,At least $5.313\%$ odd integers are a sum of one square of prime, one cube of prime and one biquadrate of prime, i.e.}
\begin{align*}
\ell=p_1^2+p_2^3+p_3^4.
\end{align*}

\section{Notation and Some Preliminary Lemmas}
\setcounter{equation}{0}
\hspace{5.8mm}For the proof of the Theorems, in this section we introduce the necessary notation and Lemmas.

Throughout this paper, by $N$ we denote a sufficiently large even integer. In addition, let $\eta < 10^{-10}$ be a fixed positive constant, and let $\varepsilon < 10^{-10}$ be an arbitrarily small positive constant not necessarily the same in different formulae.  The letter $p$, with or without subscripts, is reserved for a prime number. We use $e(\alpha)$ to denote $e^{2\pi i\alpha}$ and $e_q(\alpha) = e(\alpha/q)$. We denote by $(m,n)$ the greatest common divisor of $m$ and $n$. As usual, $\varphi(n)$ stands for  Euler's function. Let
\begin{eqnarray*}
&&P_2=\sqrt{(1-\eta) N},\,\,\,\,\,\,P_3=\left(\frac{\eta N}{2}\right)^{\frac13 },\,\,\,\,\,\,P_4=\left(\frac{\eta N}{2}\right)^{\frac{1}{4}},\,\,\,\,\,\,L=\frac{\log (\frac{N}{\log N})}{\log 2},\\
&&\mbox{}\nonumber\\
&&S_k(\alpha)=\sum_{{P_k/2} <p\leqslant P_k}(\log p) e(\alpha p^{k}),\,\,\,\,\,\, H(\alpha)=\sum_{v \leqslant L} e\left(\alpha 2^{v}\right),\,\,\,\,\,\,\mathcal{E}(\lambda) = \{\alpha \in (0,1]:|H(\alpha)| \geqslant \lambda L\}.
\end{eqnarray*}

For the application of the Hardy-Littlewood method, we need to define the Farey dissection.
For this purpose, we set
\begin{align*}
&Q_1=N^{\frac{3}{20} - 2\varepsilon},\,\,\,Q_2=N^{\frac{17}{20} + \varepsilon},
\end{align*}
and for $(a,q) = 1$, $1 \leqslant  a \leqslant  q$, put
\begin{align*}
&\mathfrak{M}(q, a)=\left(\frac{a}{q}-\frac{1}{qQ_2},\frac{a}{q}+\frac{1}{qQ_2}\right],\,\,\,\mathfrak{M}=\bigcup_{1\leqslant  q\leqslant 
Q_1}\bigcup_{{a=1}\atop{(a,q)=1}}^{q}\mathfrak{M}(q, a),\\
&\mathfrak{J}_0=\left(\frac{1}{Q_2},1+\frac{1}{Q_2}\right],\,\,\mathfrak{m}=\mathfrak{J}_0\setminus\mathfrak{M}.
\end{align*}
Then it follows from orthogonality that
\begin{align}
R(N)= & \sum_{{N = p_1^2 + p_2^2 + p_3^3 + p_4^3 + p_5^4 + p_6^4 + 2^{v_1} + 2^{v_2} + \cdots +  2^{v_k}\atop{({P_2/2})^2 <p_1,p_2 \leqslant {P_2}^2 ,\,\,({P_3/2})^3 <p_3,p_4 \leqslant {P_3}^3\atop{({P_4/2})^4 <p_5,p_6 \leqslant {P_4}^4, \,1\leqslant v_1,\cdots, v_{k}\leqslant L}}}}(\log p_1)(\log p_2)\cdots(\log p_6)\notag \\
= &\int_{0}^{1}S_2^{2}(\alpha)S_3^{2}(\alpha)S_4^{2}(\alpha) H^{k}(\alpha) e(-\alpha N ) d\alpha\notag \\
= &\left(\int_{\mathfrak{M}} + \int_{\mathfrak{m}}\right)S_2^{2}(\alpha)S_3^{2}(\alpha)S_4^{2}(\alpha)H^{k}(\alpha) e(-\alpha N ) d \alpha.
\end{align}
Now we state the lemmas required in this paper.\\

{\bf Lemma 2.1}. {\it For $(1-\eta)N \leqslant n \leqslant N$,we have
\begin{align*}
\int_{\mathfrak{M}}S_2^{2}(\alpha)S_3^{2}(\alpha)S_4^{2}(\alpha)e(-\alpha n) d\alpha = \frac{1}{2^2 \cdot 3^{2}\cdot 4^{2}} \mathfrak{S}(n) \mathfrak{J}(n)+O\left(N^{\frac{7}{6}}L^{-1}\right).
\end{align*}
Here $\mathfrak{S}(n)$ is defined as
\begin{align}
&\mathfrak{S}(n)= \sum_{q=1}^{\infty}A(n,q),\,\,\,A(n,q)=\sum_{{a=1}\atop{(a, q)=1}}^q\frac{{C_2}^{2}(q,a){C_3}^{2}(q,a){C_4}^{2}(q,a)e_{q}(-an)}{{\varphi}^{6}(q)},\\
&C_k(q, a) = \sum_{{r=1}\atop{(r, q)=1}}^qe_q(ar^k)\notag,
\end{align}
and satisfies $\mathfrak{S}(n) \gg 1$ for $n\equiv 0 \,(\bmod \,2)$. $\mathfrak{J}(n)$ is defined as
\begin{align*}
&\mathfrak{J}(n) = \sum_{{n = m_1+ m_2+ m_3+ m_4+ m_5+ m_6 \atop{({P_2/2})^2 <m_1,m_2 \leqslant {P_2}^2 ,\,\,({P_3/2})^3 <m_3,m_4 \leqslant {P_3}^3 ,\,\,({P_4/2})^4 <m_5,m_6 \leqslant {P_4}^4 }}}(m_1m_2)^{-\frac12}(m_3m_4)^{-\frac23}(m_5m_6)^{-\frac34},\notag \\
\end{align*}
and satisfies}
\begin{align*}
\mathfrak{J}(n) > (3\pi-180\eta)P_3^2P_4^2.
\end{align*}

{\bf Proof}. This follows easily from Liu [4, Lemma 2.1] and L\"{u} [6, Lemma 3.1]. \\

{\bf Lemma 2.2}. {\it For $(a,p)=1$, we have}
\begin{align*}
\mbox{i)}&\,\,|C_j(p,a)| \leqslant (j-1)p^{\frac12} + 1;\\
\mbox{ii)}&\,\,|C_3(p,a)| = -1, \mbox{for}\,\,\, p \equiv 2 \,(\operatorname{mod}3).
\end{align*}

{\bf Proof}. See Lemma 4.3 in Vaughan [7].\\

{\bf Lemma 2.3}.  {\it Let $\Xi(N,k) = \{(1-\eta)N \leqslant n \leqslant N : n = N - 2^{v_1}- 2^{v_2}- \cdots -2^{v_k},\,1\leqslant v_1,\cdots,v_k\leqslant L\}$. For $k\geqslant 16$ and $N\equiv 0 \,\,(\operatorname{mod}2)$}, {\it we have}
\begin{align*}
\mbox{i)}&\,\,\sum_{{n \in \Xi(N,k)\atop{n \equiv 0 \,(\operatorname{mod}2)}}}1 \geqslant (1-\varepsilon)L^{k};\\
\mbox{ii)}&\,\,\sum_{{n \in \Xi(N,k)\atop{n\equiv 0 \,(\operatorname{mod}2)}}}\mathfrak{S}(n) \geqslant 1.817525L^k.
\end{align*}

{\bf Proof}. For i), see Lemma 4.2 in Liu [5]. Next we give the proof of ii). It is easy to see
\begin{align}
\prod_{p\geqslant 11}\left(1+A(n, p)\right) &= \prod_{11 \leqslant p\leqslant 397}\left(1+A(n, p)\right)\prod_{397 < p\leqslant 10^6}\left(1+A(n, p)\right)\prod_{p> 10^6}\left(1+A(n, p)\right)\notag\\
&=: A_1A_2A_3.
\end{align}
For $11 \leqslant p \leqslant 397$, we directly calculate $1+A(n, p)$ by computer and obtain that
\begin{align}
A_1 \geqslant 0.916696.
\end{align}
For $397 < p\leqslant 10^6$, if $p \equiv 1 \,(\bmod \,3)$, it follows from Lemma 2.2 i) that
\begin{align}
1+A(n, p) \geqslant 1-\frac{(\sqrt{p}+1)^2(2\sqrt{p}+1)^2(3\sqrt{p}+1)^2}{(p-1)^5}.
\end{align}
Otherwise, if $p \equiv 2 \,(\bmod \,3)$, then we can deduce from Lemma 2.2 i) and ii) that
\begin{align}
1+A(n, p) \geqslant 1-\frac{(\sqrt{p}+1)^2(3\sqrt{p}+1)^2}{(p-1)^5}.
\end{align}
Consequently, on combining (2.5)-(2.6) and with the help of a computer, we deduce that
\begin{align}
A_2 &\geqslant \prod_{{397 < p\leqslant 10^6}\atop{p \equiv 1 \,(\operatorname{mod}3)}}\left(1-\frac{(\sqrt{p}+1)^2(2\sqrt{p}+1)^2(3\sqrt{p}+1)^2}{(p-1)^5}\right)\prod_{{397 < p\leqslant 10^6}\atop{p \equiv 2 \,(\operatorname{mod}3)}}\left(1-\frac{(\sqrt{p}+1)^2(3\sqrt{p}+1)^2}{(p-1)^5}\right)\notag\\
&\geqslant 0.992923.
\end{align}
For $p > 10^6$, it follows from \S 3 in Liu [5] that
\begin{align}
A_3 &\geqslant \prod_{p\geqslant 10^6}\left(1-\frac{1}{(p-1)^2}\right)^{37}\geqslant 0.999999.
\end{align}
On combining (2.3)-(2.4), (2.7)-(2.8), we obtain
\begin{align*}
\prod_{p\geqslant 11}\left(1+A(n, p)\right) \geqslant C_0:= 0.910207.
\end{align*}
Note that $A(n, p^k)=0$ for $k\geqslant 2$, $1+A(n,2)=2$ for $n\equiv 0 \,(\operatorname{mod}2)$ and $A(n, p)$ is multiplicative, we can obtain
\begin{align*}
\mathfrak{S}(n)&=\prod_{p=2}^{\infty}(1+A(n,p))\,=\prod_{2\leqslant p \leqslant 7}(1+A(n,p))\prod_{p\geqslant11}^{\infty}(1+A(n,p))\\
&\geqslant 2C_0\prod_{3\leqslant p \leqslant 7}(1+A(n,p)).
\end{align*}
For convenience, we set $q =\prod\limits_{3\leqslant p\leqslant 7}p = 105$. Then we get
\begin{align}
&\sum_{{n \in \Xi(N,k)\atop{n\equiv 0 \,(\operatorname{mod}2)}}}\mathfrak{S}(n)\geqslant 2C_0\sum_{1\leqslant j \leqslant q}\prod_{3\leqslant p\leqslant7}\left(1+A(j,p)\right)\sum_{{n \in \Xi(N,k)\atop{n\equiv 0 \,(\operatorname{mod}2)\atop{n \equiv j\,(\operatorname{mod}q)}}}}1. \label{2.3}
\end{align}
Let $S$ denote the last inner sum of $(\ref{2.3})$ and $\rho(q)$ denote the smallest positive integer $\rho$ such that $2^{\rho} \equiv 1\,(\operatorname{mod}q)$. On noting the fact that for $N\equiv 0 \,(\operatorname{mod}2)$, we have
\begin{align*}
S &= \sum_{1 \leqslant v_1,\cdots,v_k \leqslant L\atop{2^{v_1}+\cdots+2^{v_k}\equiv N\,\,(\operatorname{mod}2)\atop{2^{v_1}+\cdots+2^{v_k}\equiv N - j\,\,(\operatorname{mod}q)}}}1\,\,= \sum_{1 \leqslant v_1,\cdots,v_k \leqslant L\atop{2^{v_1}+\cdots+2^{v_k}\equiv N - j\,\,(\operatorname{mod}q)}}1\\
&= \left(\frac{L}{\rho(q)} + O(1)\right)^k\sum_{1 \leqslant v_1,\cdots,v_k \leqslant {\rho(q)}\atop{2^{v_1}+\cdots+2^{v_k}\equiv N - j\,(\operatorname{mod}q)}}1\\
&= \left(\frac L{\rho(q)}+O(1)\right)^k\frac1q\sum_{t=0}^{q-1}e\left(\frac{t(j-N)}q\right)\left(\sum_{1\leqslant s\leqslant\rho(q)}e\left(\frac{t2^s}q\right)\right)^k\\
&\geqslant\left(\frac{L}{\rho(q)}+O(1)\right)^k\frac{1}{q}\left(\rho(q)^k-(q-1)\left(\max\right)^k\right) \\
&\geqslant\frac{L^k}q\left(1-(q-1)\left(\frac{\max}{\rho(q)}\right)^k\right)+O(L^{k-1}),
\end{align*}
where
\begin{align*}
\max=\max\left\{\left|\sum_{1\leqslant s\leqslant\rho(q)}e\left(\frac{j2^s}{q}\right)\right|:1\leqslant j\leqslant q-1\right\}.
\end{align*}
 With the help of a computer, it is easy to check that
\begin{align*}
\max=6  \,\,\,\,\,\,\,\mbox{and} \,\,\,\,\,\,\, \rho(q)=12.
\end{align*}
Therefore, we can get
\begin{align}
S
&\geqslant\frac{L^k}{q}\left(1-104\times\left(\frac{1}{2}\right)^{16}\right)+O(L^{k-1})\geqslant 0.998413qL^k.
\end{align}
Moreover, we have
\begin{align}
&\sum_{1\leqslant j \leqslant q}\prod_{3\leqslant p  \leqslant 7}\left(1+A(j,p)\right)
=\,\, \prod_{3\leqslant p  \leqslant 7}\left(p + \sum_{1 \leqslant j \leqslant p}A(j,p)\right)\,\,\geqslant \,\,\prod_{3\leqslant p  \leqslant 7} p =\, q. \label{2.5}
\end{align}
On combining (2.9)-(2.11), we have
\begin{align*}
&\sum_{{n \in \Xi(N,k)\atop{n\equiv 0 \,(\operatorname{mod}2)}}}\mathfrak{S}(n) \geqslant 1.996826C_0L^k + O(L^{k-1})\geqslant 1.817525L^k. \label{2.6}
\end{align*}
Now we complete the proof of Lemma 2.3.\\

{\bf Lemma 2.4}. {\it We have}
\begin{align*}
\int_{\mathfrak{m}}|S_2(\alpha)^2S_3(\alpha)^2S_4(\alpha)^2|d\alpha\leqslant 0.514619P_3^2P_4^2.
\end{align*}

{\bf Proof}. Note that
\begin{align}
\int_{\mathfrak{m}}|S_2(\alpha)^2S_3(\alpha)^2S_4(\alpha)^2|d\alpha = \int_{0}^{1}|S_2(\alpha)^2S_3(\alpha)^2S_4(\alpha)^2|d\alpha - \int_{\mathfrak{M}}|S_2(\alpha)^2S_3(\alpha)^2S_4(\alpha)^2|d\alpha.
\end{align}
We define
\begin{align*}
\mathfrak{S}^{*}(n)&= \sum_{q=1}^{\infty}\sum_{{a=1}\atop{(a, q)=1}}^q\frac{|{C_2}^{2}(q,a){C_3}^{2}(q,a){C_4}^{2}(q,a)|}{{\varphi}^{6}(q)},
\end{align*}
and
\begin{align*}
\mathfrak{J}^{*}(n) &= \sum_{{m_1+m_2+m_3=n_1+n_2+n_3\atop{(P_2/2)^2<m_1,n_1\leqslant P_2^2\atop{(P_3/2)^3<m_2,n_2\leqslant P_3^3\atop{(P_4/2)^4<m_3,n_3\leqslant P_4^4}}}}}(m_1n_1)^{-\frac12}(m_2n_2)^{-\frac23}(m_3n_3)^{-\frac34}\notag.
\end{align*}
It follows from P.417 and Lemma 3.1 in Zhao [9] that
\begin{align}
\int_{0}^{1}|S_2(\alpha)^2S_3(\alpha)^2S_4(\alpha)^2|d\alpha \leqslant \frac{8+\eta}{2^2 \cdot 3^{2}\cdot 4^{2}}\mathfrak{S}^{*}(n)&\mathfrak{J}^{*}(n)\,\,\,\,\mbox{and}\,\,\,\,\mathfrak{S}^{*}(n)\leqslant 3.394.
\end{align}
By standard technique in analytic number theory, we can obtain
\begin{align}
\int_{\mathfrak{M}}|S_2(\alpha)^2S_3(\alpha)^2S_4(\alpha)^2|d\alpha=\frac{1}{2^2 \cdot 3^{2}\cdot 4^{2}}\mathfrak{S}^{*}(n)&\mathfrak{J}^{*}(n)+O\left(N^{\frac{7}{6}}L^{-1}\right).
\end{align}
Moreover, noting from the fact that
\begin{align*}
m_1 &= n_1 + n_2 + n_3 - m_2 -m_3\,\,\geqslant n_1- 2\eta N  \,\,\geqslant (1-4\eta)n_1,
\end{align*}
and
\begin{align*}
\sum_{(P_2/2)^2<n_1\leqslant P_2^2}n_1^{-1} &\sim 2\log 2,\,\,\,
\sum_{(P_j/2)^j<m\leqslant P_j^j}m^{1/j-1}\sim j(P_j/2)\,\, (j=3,4),
\end{align*}
we obtain
\begin{align}
\mathfrak{J}^{*}(n)&\leqslant(1-4\eta)^{-1/2}\sum_{{(P_2/2)^2<n_1\leqslant P_2^2\atop{(P_3/2)^3<m_2,n_2\leqslant P_3^3\atop{(P_4/2)^4<m_3,n_3\leqslant P_4^4}}}}(n_1)^{-1}(m_2n_2)^{-\frac23}(m_3n_3)^{-\frac34}\notag\\
&\leqslant (1+4\eta)2\log 2(3P_3/2)^2(2P_4)^2 \notag\\
&\leqslant 12.4766493P_3^2P_4^2.
\end{align}
On combining (2.12)-(2.15), we then complete the proof of Lemma 2.4.\\

{\bf Lemma 2.5}. {\it We have}
\begin{align*}
\int_{\mathfrak{m}}|S_2(\alpha)S_3(\alpha)S_4(\alpha)|^{\frac{9}{4}}d\alpha \ll N^{\frac{67}{48}+\varepsilon}.
\end{align*}

{\bf Proof}. See Lemma 4.1 in Zhang [8].\\

{\bf Lemma 2.6}. {\it We have}
\begin{align*}
\mbox{meas}\left(\mathcal{E}(\lambda)\right) \ll N_{i}^{-E(\lambda)}, \,\, \mbox{with}\,\, E(0.833783) > \frac{2}{3} + 10^{-20}.
\end{align*}

{\bf Proof}. See Lemma 5 and (3.10) in Languasco and Zaccagnini [1].

\section{Auxiliary Estimates}
\setcounter{equation}{0}
We are now equipped to establish the auxiliary estimates in this paper, and we initiate our proof by recalling the Farey dissections (2.1) that
\begin{align*}
R(N)= &\int_{0}^{1}S_2^{2}(\alpha)S_3^{2}(\alpha)S_4^{2}(\alpha)H^{k}(\alpha) e(-\alpha N ) d\alpha\notag \\
= &\left(\int_{\mathfrak{M}} + \int_{\mathfrak{m}\bigcap {\mathcal{E}{(\lambda)}}} + \int_{\mathfrak{m}\backslash {\mathcal{E}{(\lambda)}}}\right)S_2^{2}(\alpha)S_3^{2}(\alpha)S_4^{2}(\alpha) H^{k}(\alpha) e(-\alpha N ) d\alpha.
\end{align*}

{\bf Proposition 1}.  {\it  We have}
\begin{align*}
\int_{\mathfrak{M}}S_2^{2}(\alpha)S_3^{2}(\alpha)S_4^{2}(\alpha) H^{k}(\alpha) e(-\alpha N ) d\alpha  \geqslant 0.0297391P_3^2P_4^2L^{k}.
\end{align*}

{\bf Proof}. Lemmas 2.1 and 2.3 reveal that
\begin{align*}
&\int_{\mathfrak{M}}S_2^{2}(\alpha)S_3^{2}(\alpha)S_4^{2}(\alpha)H^{k}(\alpha) e(-\alpha N ) d\alpha\\
= \,\,& \sum_{n \in \Xi(N,k)}\int_{\mathfrak{M}}S_2^{2}(\alpha)S_3^{2}(\alpha)S_4^{2}(\alpha)e(-\alpha n) d\alpha\\
= \,\,& \frac{1}{2^2\cdot3^{2}\cdot4^{2}}\sum_{n \in \Xi(N,k)}\mathfrak{S}(n)\mathfrak{J}(n)  + O\left(N^{\frac{7}{6}}L^{k-1}\right)\\
\geqslant \,\,& 0.0297391P_3^2P_4^2\sum_{n \in \Xi(N,k)}1  + O\left(N^{\frac{7}{6}}L^{k-1}\right)\\
\geqslant \,\,& 0.0297391P_3^2P_4^2L^{k}.
\end{align*}

{\bf Proposition 2}.  {\it  We have}
\begin{align*}
\int_{\mathfrak{m}\bigcap {\mathcal{E}{(\lambda)}}}S_2^{2}(\alpha)S_3^{2}(\alpha)S_4^{2}(\alpha) H^{k}(\alpha) e(-\alpha N ) d\alpha \ll P_3^2P_4^2L^{k-1}.
\end{align*}

{\bf Proof}. An application of Cauchy-Schwarz inequality, Lemmas 2.5 and 2.6 yields that
\begin{align*}
&\int_{\mathfrak{m}\bigcap {\mathcal{E}{(\lambda)}}}S_2^{2}(\alpha)S_3^{2}(\alpha)S_4^{2}(\alpha) H^{k}(\alpha) e(-\alpha N ) d\alpha\\ \ll &\,\,\,L^{k}\max\limits_{\alpha \in {\mathfrak{m}}}\left(\int_{\mathfrak{m}}|S_2(\alpha)S_3(\alpha)S_4(\alpha)|^{\frac{9}{4}}d\alpha\right)^{\frac89}
\left(\int_{\mathcal{E}_{\lambda}}1d\alpha\right)^{\frac19}
\ll P_3^2P_4^2L^{k-1}.
\end{align*}
This completes the proof of Proposition 2.\\

{\bf Proposition 3}.  {\it  We have}
\begin{align*}
\int_{\mathfrak{m}\backslash {\mathcal{E}{(\lambda)}}}S_2^{2}(\alpha)S_3^{2}(\alpha)S_4^{2}(\alpha) H^{k}(\alpha) e(-\alpha N ) d\alpha \leqslant 0.514619P_3^2P_4^2\lambda^{k}L^{k}.
\end{align*}

{\bf Proof}. By Lemma 2.4, we obtain
\begin{align*}
\int_{\mathfrak{m}\backslash {\mathcal{E}{(\lambda)}}}S_2^{2}(\alpha)S_3^{2}(\alpha)S_4^{2}(\alpha) H^{k}(\alpha) e(-\alpha N ) d\alpha
&\leqslant (\lambda L)^{k}\left(\int_{\mathfrak{m}}|S_2(\alpha)S_3(\alpha)S_4(\alpha)|^2d\alpha\right)\\&\leqslant  \,0.514619P_3^2P_4^2\lambda^{k}L^{k}.
\end{align*}

\section{Proof of Theorem 1.1}
\setcounter{equation}{0}
On recalling Propositions 1,2,3, we arrive at the conclusion that
\begin{align*}
R(N) &\geqslant (0.0297391 - 0.514619\lambda^{k})P_3^2P_4^2L^{k}.
\end{align*}
When $k \geqslant 16$ and $\lambda = 0.833783$, we get
\begin{align}
R(N)>0
\end{align}
for all sufficiently large even integer $N$. Now by (4.1), the proof of Theorem 1.1 is completed.

\section{Proof of Theorem 1.2}
\setcounter{equation}{0}
For odd integers $N/8<\ell\leqslant N$, we define
\begin{align*}
r(\ell)=\sum_{{\ell=p_1^2+p_2^3+p_3^4\atop{P_2/2<p_1\leqslant P_2\atop{P_3/2<p_2\leqslant P_3\atop{P_4/2<n_2\leqslant P_4}}}}}(\log p_1)(\log p_2)(\log p_3).
\end{align*}
It follows from (2.13) and (2.15) that
\begin{align*}
\sum_{N/8<\ell\leqslant N}r^2(\ell)<\int\limits_0^1|S_2^2(\alpha)S_3^2(\alpha)S_4^2(\alpha)|d\alpha < 0.588136P_3^2P_4^2.
\end{align*}
On the one hand, we can deduce from Cauchy's inequality that
\begin{align*}
\left\{\sum_{N/8<\ell\leqslant N}r(\ell)\right\}^2&\leqslant\left\{\sum_{{N/8<\ell\leqslant N\atop{r(\ell)>0}}}1\right\}\left\{\sum_{N/8<\ell\leqslant N}r^2(\ell)\right\}\notag\\
&\leqslant 0.588136P_3^2P_4^2\left\{\sum_{{N/8<\ell\leqslant N\atop{r(\ell)>0}}}1\right\}.
\end{align*}
On the other hand, by the prime number theorem, we have
\begin{align*}
\sum_{N/8<\ell\leqslant N}r(\ell)\geqslant\sum_{P_2/2<p_1\leqslant P_2}\log p_1\sum_{P_3/2<p_2\leqslant P_3}\log p_2\sum_{P_4/2<p_3\leqslant P_4}\log p_3\geqslant\frac18(1-\varepsilon)P_2P_3P_4.
\end{align*}
Hence
\begin{align*}
\sum_{{N/8<\ell\leqslant N\atop{r(\ell)>0}}}1>\frac{1}{37.640704}(1-\varepsilon)^2P_2^2>\frac{1-\eta}{37.640704}\cdot N=\frac{1-\eta}{18.820352}\cdot\frac{N}{2},
\end{align*}
which now completes the proof of Theorem 1.2.\\
{\bf{Acknowledgement}}. The author would like to thank the anonymous referee for his/her patience and time in refereeing this manuscript.

\end{document}